\documentclass[12pt]{article}
\usepackage{amsmath, amssymb}
\begin{document}
\centerline{\bf Additional Explanatory Notes on the Analytic Proof
}\centerline{\bf of the Finite Generation of the Canonical Ring}

\bigbreak
\centerline{Yum-Tong Siu\ %
\footnote{Partially supported by a grant from the National Science
Foundation.} }

\bigbreak

\bigbreak\noindent{\sc Introduction.} This set of notes is put
together to provide some additional explanatory material on the
analytic proof of the finite generation of the canonical ring for
use by the participants in the Workshop on Minimal and Canonical
Models in Algebraic Geometry in April 16--20, 2007 at the
Mathematical Sciences Research Institute at Berkeley, California.

\bigbreak In late March 2007 before the Workshop Mihai Paun of
Strasbourg came to Harvard for two weeks to pose to me some
questions which he and other people have about the analytic proof of
the finite generation of the canonical ring which I posted in
October 2006 on the arXiv.org server [Siu 2006].  Besides orally
answering his questions I also wrote up some notes for him to give
him more precise details. This set of notes is compiled by putting
together the notes which I wrote up to answer his questions.  This
set of notes consists of two parts.

\bigbreak Part I is about how to apply the general nonvanishing
theorem to prove the precise achievement of stable vanishing order
in codimension one.  Part II gives the argument for the precise
achievement of the stable vanishing order for higher codimension.
For both Part I and Part II the new powerful tool in the analytic
proof is the use of the second case of the dichotomy of the modified
restriction of the curvature current of the metric of minimum
singularity for the canonical line bundle.

\medbreak Let $V=V_k$ be an irreducible subvariety of codimension
$k$ in the compact complex algebraic manifold $X$ of finite type
whose canonical ring is to be proved to be finitely generated such
that $V$ is an embedded branch of the stable base point set.  A {\it
modified restriction} $\Theta_V$ to $V$ of the curvature current of
the canonical line bundle is given as follows.  For some
$\eta_\nu\geq 0$ and $V_\nu$ irreducible of codimension $1$ in
$V_{\nu-1}$ ($k<\nu<0$) with $V_0=X$ and $\eta_0=0$, let $\Theta_0$
be the curvature current $\Theta_{K_X}$ on $X$ of of the metric of
minimum singularity for the canonical line bundle $K_X$ and
inductively
$$\Theta_\nu=\left(\Theta_{\nu-1}-\eta_{\nu-1}V_{\nu-1}\right)|_{V_\nu}.$$
Then $\Theta_k$ is the modified restriction $\Theta_V$ to $V$ of the
curvature current of $K_X$ for the sequence
$V_{k-1}\subset\cdots\subset V_1\subset X$ of nested subvarieties.
For the canonical decomposition
$$
\Theta_V=\sum_{j=1}^J\gamma_j\left[Y_j\right]+R
$$
of $\Theta_V$ on $V$ there are two cases. The first case is either
$J=\infty$ or $R\not=0$ and the second case is $J<\infty$ and $R=0$,
where $Y_j$ is a subvariety of codimension one in $V$ and $R$ is a
closed positive current on $V$ whose Lelong number is zero outside a
countable union of subvarieties of codimension at least two in $V$.
For Part I on how to apply the general nonvanishing theorem to prove
the precise achievement of stable vanishing order in codimension
one, the result of the use of this new powerful tool is that the
second case of the dichotomy must always eventually occur at some
positive dimensional $V$, because if the first case occurs all the
time down to dimension zero then there is some improvement in the
stable vanishing at some point of $V$ in the multi-directions
defined by the sequence of nested subvarieties.  The second case of
the dichotomy gives a explicitly constructed section on $V$ (unique
up to a nonzero constant factor) which belongs to the multiplier
ideal sheaf corresponding to $\Theta_V$.

\medbreak In Part I of these additional explanatory notes, in order
to make the arguments of the general nonvanishing theorem more
transparent, I treat separately the two cases of the dichotomy. We
use the process of using the techniques for the Fuijta conjedture by
constructing singular metrics successively to cut down on the
dimension of (the projection of) the zero-set of the multiplier
ideal sheaf until we end up with the inevitable second case of the
dichotomy and in the second case of the dichotomy we use the
sections explicitly constructed from the special form of the
canonical decomposition of the modified restriction of the curvature
current. In extending the explicitly constructed sections all the
way back to the ambient manifold $X$ we introduce the technique of
{\it constrained} minimum center of log canonical singularity.

\medbreak In Part II of these additional explanatory notes about the
argument for the precise achievement of the stable vanishing order
for higher codimension, we give its illustration in the low
dimensional cases of complex surfaces and complex threefolds.  The
illustration in the low dimensional cases avoids the encumberment of
complicated notations in the case of general dimension and at the
same time contains the essence of the argument of the general case.

\medbreak For the case of a threefold $X$ in Part II, the key case
of higher codimension is that of an irreducible curve $C$ which is
an embedded curve in the stable base point set. The important
ingredient of using the second case of the dichotomy is the
following.  For simplicity assume that there is no codimension-one
base point set.

\medbreak Suppose $s_1,\cdots,s_k$ are pluricanonical sections on
$X$ such that $C$ is cut out by the family $S_\sigma$ of surfaces
defined by pluricanonical sections $s_\sigma=\sum_{j=1}^k\sigma_j
s_j$ parametrized by
$\sigma=\left(\sigma_1,\cdots,\sigma_k\right)\in{\mathbb
C}^k-\left\{0\right\}$. Suppose $s$ is a pluricanonical section
whose restriction to some $S_\sigma$ achieves the {\it restriction}
of the stable vanishing order $\alpha_\sigma$ on $S_\sigma$ at some
point $P_\sigma$ of $C$. Then, after taking away the restriction of
the stable vanishing order $\alpha_\sigma$ of $s|_{S_\sigma}$ across
$C$, the resulting section
$\frac{s|_{S_\sigma}}{\left(s_{C,\sigma}\right)^{\alpha_\sigma}}$ on
$S_\sigma$ (where $s_{C,\sigma}$ is the canonical section of $C$ in
$S_\sigma$), when restricted to $C$, is the section (up to a nonzero
constant factor) explicitly constructed from the second case of the
dichotomy for $C$ in $S_\sigma$. This means that at every point of
$C$ outside the finite zero-set $Z_\sigma$ of the resulting section
$\frac{s|_{S_\sigma}}{\left(s_{C,\sigma}\right)^{\alpha_\sigma}}$ on
$C$, the section
$\frac{s|_{S_\sigma}}{\left(s_{C,\sigma}\right)^{\alpha_\sigma}}$
achieves the restriction of the stable vanishing order
$\alpha_\sigma$ on $S_\sigma$.

\medbreak Suppose for every parameter $\sigma$ there is a parameter
$\tau$ and some point $P_{\sigma,\tau}$ in $C$ such that the
restriction of the pluricanonical section $s_\tau$ to $S_\sigma$
achieves the {\it restriction} of the stable vanishing order on
$S_\sigma$ at the point $P_{\sigma,\tau}$ of $C$.  Let $Z$ be the
finite subset of $C$ such that as a point $P$ varies in $C$ outside
$Z$ the Artinian subscheme defined by $s_1,\cdots,s_k$ in the normal
direction of $C$ at $P$ does not change.  In other words, every
point $Q\in C-Z$ admits an open neighborhood $U$ in $C-Z$ and a
holomorphic family of local nonsingular complex surfaces $M_P$ in
$X$ parametrized by $P\in U$ such that each $M_P$ intersecting $C$
normally at $P$ and
$$
\dim_{\mathbb C}{\mathcal O}_{M_P}\left/\sum_{j=1}^k {\mathcal
O}_{M_P}\left(s_j|_{M_P}\right)\right.
$$
is independent of $P\in U$.  From the above discussion about $s$
being the explicitly constructed section and about the achievement
of the restriction to $S_\sigma$ of stable vanishing order at points
of $C$ outside $Z_\sigma$ it follows that the stable vanishing order
is achieved at points of $C-Z$.

\medbreak This use of the sections, explicitly constructed from the
curvature current in the second case of the dichotomy, in order to
come up with a finite subset $Z$ of $C$ can be regarded as the
higher-codimensional way of constructing sections of {\it vector
bundles} explicitly from the second case of the dichotomy of the
modified restriction of the curvature current.  Note that the
codimension one case for a surface is invoked for us to get the
rationality of the restriction of the stable vanishing order to a
surface and its achievement on the surface at some point of $C$.
This part of the hard work already done for the codimension one case
is used here.  This {\it vector bundle} section on $C$ explicitly
constructed from the second case of the dichotomy of the modified
restriction of the curvature current is obtained, so to speak, by
piecing together the {\it line bundle} sections on $C$ explicitly
constructed by using the surface slices $S$ from the second case of
the dichotomy of the modified restriction to the surface slice $S$
of the curvature current.  This is the reason why the usual worry
about interminable blowups to transform into hypersurfaces newly
appearing higher-codimensional base point sets is not a problem
here.

\bigbreak Before presenting Part I and Part II, we would like to
make two remarks.  The first remark is about the finite generation
of the canonical ring without the assumption of finite type.  In
analysis, when $K_X$ is not big, adding a line bundle $E$ to $K_X$
to have a big bundle $K_X+E$ and working with $K_X+E$ instead of
$K_X$ would only mean a standard modification of the argument for
the finite generation of the canonical ring for the case of general
type. (This is actually used in the form of $K_X-\gamma Y$, with
$\gamma$ being the stable vanishing order across the hypersurface
$Y$, in the proof of the finite generation of the canonical ring for
the case of general type.) However, getting rid of $E$ afterwards
needs the justification of a rather involved limiting process.  Such
a limiting process is of the same nature as the limiting process
which is used to extend the proof of the deformational invariance of
the plurigenera for the case of general type [Siu 1998] to the
general algebraic case [Siu 2002].

\medbreak For the finite generation of the canonical ring, when
$K_X$ is big, we introduce the curvature $\Theta_{K_X}$ of the
metric of minimum singularity constructed from pluricanonical
sections. When $K_X$ is not big, we use an ample line bundle $A$ on
$X$ and the limit curvature current
$$
\tilde\Theta_{K_X}:=\lim_{m\to\infty}\Theta_{K_X+\frac{1}{m}A}
$$
to replace $\Theta_{K_X}$.  All the diophantine approximation
arguments deal with the modified restriction of the canonical
decomposition of $\tilde\Theta_{K_X}$. This makes it possible to go
ahead with the limiting process.

\medbreak In the case of the deformational invariance of the
plurigenera the limiting process is handled by abandoning the metric
of minimum singularity and replacing it by another metric with
maximum allowable singularity which is still good enough for the
finiteness of the $L^2$ norm of the section to be extended, because
there is a need to bound the dimension of the space of $L^2$
holomorphic $m$-canonical sections with the bound independent of $m$
or at least growing more slowly than $m$. It is not yet clear what
the analogous procedure is for the problem of the finite generation
of the canonical ring.

\medbreak The second remark is to answer the question why one does
not simply take a pluricanonical section and get the finite
generation of the canonical ring of its divisor and then use
extension of pluricanonical sections from that divisor to the entire
ambient space.  This procedure is explained in \S8 of the posted
notes of my analytic proof [Siu 2006].  The extension technique for
pluricanonical sections in the problem of deformational invariance
of the plurigenera was developed with that particular application in
mind [Siu 1998, Siu 2006]. Though the pluricanonical extension
technique for deformational invariance of the plurigenera is for
extension from the fiber at a point to the open unit disk, in
analysis it does not matter whether one deals with a compact
manifold or a Stein manifold if the Stein manifold is a Zariski open
subset and one keeps track of the $L^2$ bounds.

\medbreak However, the situation of the singularities in the
pluricanonical divisor poses technical problems in analysis.  In my
article in the proceedings of a 2001 conference in Hanoi [Siu 2004]
I addressed the problem of pluricanonical extension from a singular
divisor.  There are technical problems concerning finite $L^2$
bounds.  The difficulty in handling such technical problems
concerning finite $L^2$ bounds is the same as getting the
comparability of the metric of minimum singularity for $K_X$ and
that from an appropriate truncated finite sum in the application of
Skoda's theorem of ideal generation [Skoda 1972]. Moreover, the
Ohsawa-Takegoshi-type extension theorem used in the pluricanonical
extension is closely related to Skoda's theorem of ideal generation
(see, for example, Ohsawa's article in the Festschrift for Grauert's
70th birthday [Ohsawa 2002]) and the constant in Skoda's theorem is
more precise. That is why we use the approach by Skoda's theorem
instead of the essentially equivalent approach of pluricanonical
extension.

\medbreak It seems that at this point any approach to the problem of
finite generation of the canonical ring has to rely, as starting
point, on the technique of pluricanonical extension. Any difference
in the different approaches are more a matter of technical handling
of the singularities of the pluricanonical divisor.

\medbreak In the definition of a stable metric in [Siu 2006, (6.1)]
there are some typographical errors and inaccuracies. The correct
definition is as follows.  A metric $e^{-\kappa}$ of a line bundle
$L$ on a compact K\"ahler manifold $M$ is said to be {\it stable} if
$\kappa$ is locally plurisubharmonic and there exists some
$\varepsilon>0$ with the following property. If $U$ is an open
neighborhood of a point $P\in M$, and $\varphi$ and $\psi$ are
plurisubharmonic functions on $U$ such that the total mass, with
respect to the K\"ahler form of $M$, of the sum of the two closed
positive $(1,1)$-currents $\Theta_\varphi$ and $\Theta_\psi$ is less
than $\epsilon$ and if $\kappa-\varphi$ is plurisubharmonic, then
there exists an open neighborhood $U^\prime$ of $P$ in $M$ such that
the multiplier ideal sheaf ${\mathcal I}_{\kappa+\psi}$ of the
metric $e^{-\kappa-\psi}$ on $U^\prime$ is equal to the multiplier
ideal sheaf ${\mathcal I}_{\kappa-\varphi}$ of the metric
$e^{-\kappa+\varphi}$ on $U^\prime$.

\eject\centerline{\bf PART I} \centerline{\bf How to Apply
Nonvanishing Theorem to Precisely Achieve}\centerline{\bf Stable
Vanishing Order in Codimension One}

\bigbreak

\bigbreak The proof of the precise achievement of stable vanishing
order in codimension one has the following ingredients.
\begin{itemize}\item[(a)] The techniques for the Fujita conjecture
which consists of
\begin{itemize}\item[(i)] constructing singular metrics with curvature current of strict
positive lower bound whose multiplier ideal sheaf has high vanishing
order at a prescribed point,
\item[(ii)] blowing up the zero-set of the multiplier ideal sheaf of the
new singular metric,
\item[(iii)] repeating the procedure so that the zero-set of the
multiplier ideal sheaf of the new singular metric on a blow-up space
is projected down to lower and lower dimension in the original
manifold, until one gets to the case of a singular point in the
original manifold, and
\item[(iv)] extending a section defined on the single point to over
all of the original manifold by the vanishing theorem of
Kawamata-Viehweg-Nadel.
\end{itemize}
Because of the need of adding one canonical line bundle in the
application of the vanishing theorem of Kawamata-Viehweg-Nadel, the
singular metric is for a multiple of the line bundle in question
minus the canonical line bundle.
\item[(b)] Introducing a dichotomy depending on the canonical
decomposition of some curvature current so that
\begin{itemize}\item[(i)] in the first case of the dichotomy the
technique of the Fujita conjecture for the construction of singular
metric can be carried out, and
\item[(ii)] in the second case of the dichotomy a section can be
explicitly constructed which can be extended to all of the original
manifold.
\end{itemize}
Once we get to the second case of the dichotomy there is no need to
use the technique of the Fujita conjecture to construct any more
singular metric.  The process is stopped and complete.  For the
precise achievement of stable vanishing order, the second case of
the dichotomy must arise, otherwise the stable vanishing order can
be improved, which contradicts the definition of a stable vanishing
order.  The explicitly constructed section is rather rigid in the
sense that there is no choice and the section comes from the
curvature current in a rather unique way.  This uniqueness and
rigidity of the explicitly constructed section will be a key
ingredient in the proof of the precise achievement of the stable
vanishing order in higher codimension.
\item[(c)] The new technique of constrained minimum center of log
canonical singularity, whose motivation and precise application will
be explained below.  This new technique is needed, because of the
undesirable additional vanishing order in the process of
constructing singular metrics, which will be explained in detail
below. \item[(d)] Use of diophantine approximation to handle
irrational coefficients.  Since the use of diophantine approximation
has already been described in detail in the posted notes of my proof
of the finite generation of the canonical ring [Siu 2006], we will
not discuss diophantine approximation in these additional
explanatory notes and we just assume that the relevant coefficients
are known to be rational numbers.
\end{itemize}
As the first step we will consider the first case of the dichotomy
and describe how to construct an appropriate singular metric.  Then
we will give the motivation for constrained minimum center of log
canonical singularity by reviewing the goal and the strategy of the
standard technique for the Fujita conjecture. Finally we consider
the second case of the dichotomy where a rather rigid section is
explicitly constructed.

\bigbreak\noindent{\sc Proposition}{\it (Construction of Metric with
Multiplier Ideal Sheaf Vanishing to High Order at a Prescribed Point
for the First Case of the Dichotomy of the Curvature Current).} Let
$M$ be a compact complex projective algebraic manifold of complex
dimension $n$. Let $L$ be a holomorphic line bundle on $M$ with a
(possibly singular) metric $e^{-\varphi}$ along its fibers whose
curvature current $\Theta_\varphi$ is a closed positive
$(1,1)$-current. Let
$$
\Theta_\varphi=\sum_{j=1}^J\tau_j\left[V_j\right]+R
$$
be the canonical decomposition of $\Theta_\varphi$, where
$J\in{\mathbb N}\cup\left\{0,\infty\right\}$ and the Lelong number
of $R$ is zero outside a countable union $Z$ of subvarieties of
codimension at least two in $M$ and $V_j$ is an irreducible
hypersurface in $M$ and $\tau_j>0$.  Assume that either $J=\infty$
or $R\not=0$ (that is, one is in the first case of the dichotomy for
the curvature current $\Theta_\varphi$). Assume that for some
positive integer $p_0$ there is a (possibly singular) metric
$e^{-\chi}$ along the fibers of $p_0L-K_M$ which is stable and whose
curvature current $\Theta_\chi$ is a closed positive $(1,1)$-current
which dominates some strictly positive smooth $(1,1)$-form on $M$.
Let $P_0$ be a point of $M$ such that the Lelong number of
$\Theta_\varphi$ is zero at $P_0$. Let $q\in{\mathbb N}$.  Then for
some sufficiently divisible $m\in{\mathbb N}$ the line bundle
$\left(m+p_0\right)L-K_M$ admits a metric $e^{-\tilde\chi}$ whose
curvature current dominates some strictly positive smooth
$(1,1)$-form on $M$ such that its multiplier ideal sheaf ${\mathcal
I}_{\tilde\chi}$ is contained in the maximum ideal of $M$ at $P_0$
raised to the $q$-th power. Moreover, if $M$ is a hypersurface in
some compact complex algebraic manifold $X$ of general type so that
$L$ is the restriction of some line bundle $\tilde L$ on $X$ and
$K_M=K_X+M$ and the metric $e^{-\chi}$ is defined by a convergent
infinite sum of multi-valued holomorphic sections of $p_0\tilde
L-K_X-M$ over $X$, then the metric $e^{-\tilde\chi}$ can be chosen
to be defined also by a convergent infinite sum of multi-valued
holomorphic sections of $\left(p+p_0\right)\tilde L-K_X-M$ over $X$.

\bigbreak\noindent{\sc Proof.} The idea of the proof is to use the
techniques for Fujita's conjecture (see, for example, [Angehrn-Siu
1995]).

\bigbreak\noindent{\it Slicing by an Ample Divisor.} Let $A$ be a
very ample line bundle over $M$ such that $A-K_M$ is ample. Let
$h_A$ be a smooth metric of $A$ whose curvature form $\omega_A$ is
positive on $M$.  We assume that $A$ is chosen to be sufficiently
ample so that for each point $P\in M$ the proper transform of $A$ in
the manifold obtained from $M$ by blowing up $P$ is still very
ample. This technical assumption will enable us to choose a generic
element of $\Gamma\left(M,\,A\right)$ vanishing at $P_0$ which is
not a zero-divisor of a prescribed coherent ideal sheaf.

\medbreak Let $p$ and $k$ be positive integers and we will impose
more conditions on $p$ and $k$ later. Let $s_1$ be a generic element
of $\Gamma\left(M,\,A\right)$ vanishing at $P_0$ so that the short
exact sequence
$$
\displaylines{\qquad 0\to {\mathcal
I}_{p\varphi+\chi}\left(\left(p+p_0\right)L-K_M
+kA\right)\stackrel{\theta_{s_1}}{\longrightarrow}{\mathcal
I}_{p\varphi+\chi}\left(\left(p+p_0\right)L-K_M+\left(k+1\right)A\right)\hfill\cr\hfill\to\left({\mathcal
I}_{p\varphi+\chi}\left/s_1{\mathcal
I}_{p\varphi+\chi}\right.\right)\left(\left(p+p_0\right)L-K_M+\left(k+1\right)A\right)\to
0\qquad\cr}
$$
is exact, where $\theta_{s_1}$ is defined by multiplication by
$s_1$.  Let $M_1$ be the zero-set of $s_1$ and
$$
{\mathcal O}_{M_1}=\left({\mathcal O}_M\left/s_1{\mathcal
O}_M\right.\right)|_{M_1},
$$
which we can assume to be regular with ideal sheaf equal to
$s_1{\mathcal O}_M$ because $s_1$ is generic element of
$\Gamma\left(M,A\right)$ vanishing at $P_0$.  By choosing $s_1$
generically we can also assume that ${\mathcal
I}_{\left(p\varphi+\chi\right)|_{M_1}}={\mathcal
I}_{p\varphi+\chi}\left/s_1{\mathcal I}_{p\varphi+\chi}\right.$ and
${\mathcal I}_{\left(p\varphi\right)|_{M_1}}={\mathcal
I}_{p\varphi}\left/s_1{\mathcal I}_{p\varphi}\right.$ and ${\mathcal
I}_{\chi|_{M_1}}={\mathcal I}_\chi\left/s_1{\mathcal
I}_\chi\right.$. We use $\chi\left(\cdot,\,\cdot\right)$ to denote
the arithmetic genus which means
$$
\chi\left(\cdot,\,\cdot\right)=\sum_{\nu=0}^\infty(-1)^\nu\dim_{\mathbb
C}H^\nu\left(\cdot,\,\cdot\right).
$$
From the long cohomology exact sequence of the above short exact
sequence we obtain
$$
\displaylines{\chi\left(M,\,{\mathcal
I}_{p\varphi+\chi}\left(\left(p+p_0\right)L-K_M+\left(k+1\right)A\right)\right)=\cr
\chi\left(M,\,{\mathcal
I}_{p\varphi+\chi}\left(\left(p+p_0\right)L-K_M+kA\right)\right)+\chi\left(M_1,\,{\mathcal
I}_{\left(p\varphi+\chi\right)|_{M_1}}\left(\left(p+p_0\right)L-K_M
+\left(k+1\right)A\right)|_{M_1}\right).\cr}
$$
Since $A-K_M$ is ample and $2A-K_{M_1}=A-K_M$ is also ample, when we
assume $k\geq 1$, by the theorem of Kawamata-Viehweg-Nadel
$$
\displaylines{H^\nu\left(M,{\mathcal
I}_{p\varphi+\chi}\left(\left(p+p_0\right)L-K_M+kA\right)\right)=0\quad{\rm
for\ }\nu\geq 1,\cr H^\nu\left(M_1,{\mathcal
I}_{\left(p\varphi+\chi\right)|_{M_1}}\left(\left(\left(p+p_0\right)L-K_M
+\left(k+1\right)A\right)|_{M_1}\right)\right)=0\quad{\rm for\
}\nu\geq 1\cr}
$$
so that
$$
\displaylines{\dim_{\mathbb C}\Gamma\left(M,\,{\mathcal
I}_{p\varphi+\chi}\left(\left(p+p_0\right)L-K_M+\left(k+1\right)A\right)\right)\cr
\qquad\qquad\qquad=\dim_{\mathbb C}\Gamma\left(M,\,{\mathcal
I}_{p\varphi+\chi}\left(\left(p+p_0\right)L-K_M+kA\right)\right)\hfill\cr\hfill+\dim_{\mathbb
C}\Gamma\left(M_1,\,{\mathcal
I}_{\left(p\varphi+\chi\right)|_{M_1}}\left(\left(\left(p+p_0\right)L-K_M+\left(k+1\right)A\right)|_{M_1}\right)\right)\cr
\geq\dim_{\mathbb C}\Gamma\left(M_1,\,{\mathcal
I}_{\left(p\varphi+\chi\right)|_{M_1}}\left(\left(\left(p+p_0\right)L-K_M+\left(k+1\right)A\right)|_{M_1}\right)\right)
.\cr }
$$

\bigbreak\noindent{\it Slicing by Ample Divisors Down to a Curve.}
Instead of one single element $s\in\Gamma\left(M,A\right)$, we can
choose generically
$$
s_1, \cdots, s_{n-1} \in\Gamma\left(M,\,A\right)
$$
all vanishing at $P_0$ so that inductively for $1\leq\nu\leq n-1$
the common zero-set $M_\nu$ of $s_1,\cdots,s_\nu$ with the structure
sheaf
$$
{\mathcal O}_{M_\nu}:=\left({\mathcal O}_M\left/\sum_{j=1}^\nu
s_j{\mathcal O}_M\right)\right|_{M_\nu}
$$
is regular and we end up with the inequality
$$
\displaylines{\dim_{\mathbb C}\Gamma\left(M,\,{\mathcal
I}_{p\varphi+\chi}\left(\left(p+p_0\right)L-K_M+\left(k+n-1\right)A\right)\right)\cr\geq
\dim_{\mathbb C}\Gamma\left(M_{n-1},\,{\mathcal
I}_{\left(p\varphi+\chi\right)|M_{n-1}}\left(\left(\left(p+p_0\right)L-K_M+\left(k+n-1\right)A\right)|_{M_{n-1}}\right)\right).\cr}
$$
Since $M_{n-1}$ is a curve, all coherent ideal sheaves on it are
principal and are locally free and they come from holomorphic line
bundles.

\medbreak We would like to remark also that this particular step of
slicing by $n-1$ ample divisors to get down to a curve roughly
corresponds to the step in Shokurov's proof of his non-vanishing
theorem [Shokurov 1985] where he takes the product of his
numerically effective divisor in his $n$-dimensional manifold with
the $(n-1)$-th power of a numerically effective big line bundle.

\bigbreak\noindent{\it Application of the Theorem of Riemann-Roch to
a Curve and Comparing Contributions from the Curvature Current and
the Multiplier Ideal Sheaves.}  Let $b$ be the Chern class of the
line bundle on $M_{n-1}$ defined by the multiplier ideal sheaf
${\mathcal I}_{\chi|_{M_{n-1}}}$ of the restriction to $M_{n-1}$ of
the metric $e^{-\chi}$. Let $c$ be the nonnegative number
$$\int_{M_{n-1}}R=\int_M
R\wedge\left(\omega_A\right)^{n-1}.$$ Then
$$
\displaylines{(\dagger)\qquad\qquad\dim_{\mathbb
C}\Gamma\left(M,\,{\mathcal
I}_{p\varphi+\chi}\left(\left(p+p_0\right)L-K_M+\left(k+n-1\right)A\right)\right)\hfill\cr\geq
\dim_{\mathbb C}\Gamma\left(M_{n-1},\,{\mathcal
I}_{\left(\left.p\varphi+\chi\right|M_{n-1}\right)}\left(\left(\left(p+p_0\right)L-K_M+\left(k+n-1\right)A\right)|_{M_{n-1}}\right)\right)\cr
\geq 1-{\rm genus}\left(M_{n-1}\right)+b+
\left(k+n-1\right)A^{n-1}M_{n-1}
\cr+\sum_{j=1}^J\left(p\tau_j-\left\lfloor
p\tau_j\right\rfloor\right) V_j\cdot A^{n-1}+p\int_{M_{n-1}}R,\cr}
$$
where the last identity is from the theorem of Riemann-Roch applied
to the regular curve $M_{n-1}$ and the locally free sheaf
$${\mathcal
I}_{\left(\left.p\varphi+\chi\right|M_{n-1}\right)}\left(\left(\left(p+p_0\right)L-K_M+\left(k+n-1\right)A\right)|_{M_{n-1}}\right)$$
on $M_{n-1}$. From the assumption that $J=\infty$ or $R\not=0$, we
conclude that the right-hand side of $(\dagger)$ goes to $\infty$ as
$p$ goes to $\infty$ through an appropriate sequence, where for the
case of $J=\infty$ and $R=0$ a diophantine argument has to be used
whereas for the case $R\not=0$ we simply need to use $c>0$.

\bigbreak\noindent{\it Construction of Sections with Extra Vanishing
Order from Dimension Counting and Construction of Metrics by
Canceling Contributions from Ample Divisors by Using the General
Type Property.} For any $\ell\in{\mathbb N}$ the number of terms in
a polynomial of degree $\ell$ in $d$ variables is ${d+\ell\choose
d}$. Take a positive integer $N$ and we will impose more condition
on $N$ later.  By the behavior of the right-hand side of $(\dagger)$
as $p\to\infty$, there exists $p\in{\mathbb Z}$ such that
$$
\dim_{\mathbb C}\Gamma\left(M, \,{\mathcal
I}_{p\varphi+\chi}\left(\left(p+p_0\right)L-K_M+\left(k+n-1\right)A\right)\right)\geq
1+{n+N\left(k+n-1\right)q\choose n}
$$
and we can find some non identically zero element $s$ of
$$\Gamma\left(M, \,{\mathcal
I}_{p\varphi+\chi}\left(\left(p+p_0\right)L-K_M+\left(k+n-1\right)A\right)\right)$$
which vanishes to order at least $N\left(k+n-1\right)q$ at $P_0$ so
that $s^{\frac{1}{N\left(k+n-1\right)}}$ is a multi-valued
holomorphic section of the ${\mathbb Q}$-line-bundle
$\frac{p}{N\left(k+n-1\right)}\,L+\frac{1}{N}\,A$ over $M$ which
vanishes to order at least $q$ at $P_0$.  We assume that $N$ is
chosen so large that the curvature current $\Theta_\chi$ dominates
$\frac{2}{N}\omega_A$.  Let $\hat p$ to be the round-up of
$\frac{p}{N\left(k+n-1\right)}$ and $\delta_p=\hat
p-\frac{p}{N\left(k+n-1\right)}$.  We introduce the metric
$$
e^{-\tilde\chi}:=\frac{e^{-\chi-\delta_p\varphi}}{\left(h_A\right)^{\frac{1}{N}}\left|s\right|^{\frac{2}{Nk}}}
$$
of $\left(p+p_0\right)L-K_M$ so that the multiplier ideal of
${\mathcal I}_{\tilde\chi}$ at $P_0$ is contained in
$\left({\mathfrak
m}_{M,P_0}\right)^{{}^{\left\lfloor\frac{q}{n}\right\rfloor}}$.
Q.E.D.

\bigbreak\noindent{\sc Remark} {\it on Application of the
Proposition on Construction of Metric for the First Case of the
Dichotomy.}  The above application is applied in the following
manner.  Let $X$ be the the compact complex algebraic manifold of
finite type whose canonical ring is to be proved to be finitely
generated.  Let $Y$ be a hypersurface in $X$ across which the stable
vanishing order is $\gamma>0$.  Let $P_0$ be a generic point of $Y$.
We start out with $M=Y$ and $L=K_X-\gamma Y$. After we get the new
metric $e^{-\tilde\chi}$, we use an interpolation between
$e^{-\tilde\chi}$ and $e^{-p\varphi-\chi}$ and a slight modification
to get to a minimum center of log canonical singularity which, after
blow-up, projects down to a proper subvariety of $M$ containing
$P_0$.  Then we replace $X$ by its blow-up and replace $M$ by the
new minimum center of log canonical singularity and replace $L$ by
its pullback to the blowup of $X$.  We continue doing this until we
inevitably come to the second case of dichotomy eventually as
explained in the Introduction.

\bigbreak\noindent{\sc Remark} {\it on the Second Case of the
Dichotomy.}  Suppose
$$
\Theta_\varphi=\sum_{j=1}^J\tau_j\left[V_j\right]
$$
with $J<\infty$.  Then we can explicitly construct a section $s_0$
of $pL$ over $M$.  The reason why a minimal center of log canonical
singularity is used in the techniques for the Fujita conjecture is
to make sure that when we take the subspace defined by the
multiplier ideal sheaf, the subspace has a reduced structure.  In
our case we have to introduce the notion of {\it constrained}
minimal center of log canonical singularity so that the center is
not completely contained in the zero-set of $s_0$.  For that we have
to pay the price that the subspace defined by the multiplier ideal
sheaf may have an unreduced structure, but the set where nonzero
nilpotent elements of its structure sheaf occurs is contained in the
zero-set of $s_0$.  By raising $s_0$ to a sufficiently high power,
we can handle the unreduced structure and get the extension of a
sufficiently high power of $s_0$ to the ambient manifold $X$ by
using the vanishing theorem from the metric of $pL-K_M$.  We are
going to elaborate on this by reviewing the goal of the techniques
for the Fujita conjecture and the use of minimal center of log
canonical singularity and also how we are naturally and by necessity
led to the concept of constrained minimal center of log canonical
singularity.

\bigbreak\noindent{\it Main Idea of the Techniques for the Fujita
Conjecture.}  For the discussion about the main idea of the
technique for the Fujita conjecture, we forget the above meaning of
$X$ and $L$ and for the time being use the symbols $X$ and $L$ in
another context.  We will so indicate when we later return to the
above meaning of $X$ and $L$.  The goal of the technique for the
Fujita conjecture is to find global sections to generate some
positive power $mL$ of a line bundle $L$ over a compact complex
algebraic manifold $X$.

\medbreak For a more general setting, the goal is to find global
sections to globally generate ${\mathcal J}(mL)$ over $X$ for some
given coherent ideal sheaf ${\mathcal J}$ at points outside some
given subvariety $Z$ of $X$. The problem of proving the finite
generation of the canonical ring by verifying the precise
achievement of stable vanishing orders actually involves this more
general setting.  There the even more complicated situation of
supremum norm is used instead of just the $L^2$ norm.  However, for
the sake of simplicity in our discussion of the main idea of the
techniques of the Fujita conjecture, we stick with the simpler goal
of find global sections to generate some positive power $mL$ of a
line bundle $L$ over a compact complex algebraic manifold $X$.

\medbreak The main idea of the technique is to focus on the
subvariety where the global generation fails.  We take a basis of
$s_1,\cdots,s_k\in\Gamma\left(X,mL\right)$ and let $Y$ be their
common zero-set so that global generation precisely fails at points
of $Y$.  The main idea of the technique is simply to focus on $Y$ if
$Y$ is nonempty.  We seek elements of $\Gamma\left(Y,mL|_Y\right)$
which are not identically zero and then extend them to elements of
$\Gamma\left(X,mL\right)$, which would then contradict the
definition of $Y$.  Hopefully the extension of elements of
$\Gamma\left(Y,mL|_Y\right)$ to elements of
$\Gamma\left(X,mL\right)$ could be done by the vanishing of some
appropriate first cohomology group.  Usually this first cohomology
comes from the vanishing theorem of Kawamata-Viehweg-Nadel.  We seek
a metric $e^{-\varphi}$ of $mL-K_X$ so that\begin{itemize}\item[(i)]
the zero-set of its multiplier ideal sheaf ${\mathcal I}_\varphi$ is
$Y$, and \item[(ii)] the curvature current of $e^{-\varphi}$
dominates some strictly positive smooth $(1,1)$-form on
$X$.\end{itemize} In such a case we have
$$
H^1\left(X,{\mathcal I}_\varphi\left(mL\right)\right)=0
$$
from the vanishing theorem of Kawamata-Viehweg-Nadel and the map
$$
\Gamma\left(X,mL\right)\to\Gamma\left(Y,\left({\mathcal
O}_X\left/{\mathcal I}_\varphi\right.\right)\left(mL\right)\right)
$$
is surjective.  The next step is to come up with some element of
$$\Gamma\left(Y,\left({\mathcal O}_X\left/{\mathcal
I}_\varphi\right.\right)\left(mL\right)\right)$$ which induces a non
identically zero element of $\Gamma\left(Y,mL|_Y\right)$.

\medbreak It is at this point that the question of a possibly {\it
unreduced} complex structure ${\mathcal O}_X\left/{\mathcal
I}_\varphi\right.$ arises.  It means that the structure sheaf
${\mathcal O}_X\left/{\mathcal I}_\varphi\right.$ for $Y$ may have
nonzero nilpotent elements. This is the case when ${\mathcal
I}_\varphi$ is a {\it proper} subsheaf of the {\it full} ideal sheaf
${\mathcal I}_Y$ of $Y$ and is not equal to ${\mathcal I}_Y$.  An
element of ${\mathcal I}_Y$ which is not in ${\mathcal I}_\varphi$
would yield a nonzero nilpotent element of the structure sheaf
${\mathcal O}_X\left/{\mathcal I}_\varphi\right.$ for $Y$.  When we
have an unreduced structure ${\mathcal O}_X\left/{\mathcal
I}_\varphi\right.$ for $Y$, it is more difficult to produce some
element of
$$\Gamma\left(Y,\left({\mathcal O}_X\left/{\mathcal
I}_\varphi\right.\right)\left(mL\right)\right)$$ which induces a non
identically zero element of $\Gamma\left(Y,mL|_Y\right)$.  To handle
the problem of unreduced structure sheaf, the technique is to use
minimum centers of log canonical singularity whose role we are going
to explain.

\bigbreak\noindent{\it Minimum Center of Log Canonical Singularity.}
The idea is to seek a metric $e^{-\psi}$ for $mL-K$ which is less
singular than $e^{-\varphi}$ so that the multiplier ideal sheaf
${\mathcal I}_\psi$ of $e^{-\psi}$ contains the multiplier ideal
sheaf ${\mathcal I}_\varphi$ of $e^{-\varphi}$.  This procedure
usually involves the interpolation of two metrics and a slight
modification of the result of the interpolation.  We seek to make
the metric of $e^{-\psi}$ for $mL-K$ to be as least singular as
possible, with just enough singularity to make the multiplier ideal
sheaf ${\mathcal I}_\psi$ of $e^{-\psi}$ not equal to ${\mathcal
O}_X$.  Let $Y^\prime$ be the support of ${\mathcal
O}_X\left/{\mathcal I}_\psi\right.$.  This kind of least or minimum
singularity for the choice of $e^{-\psi}$ gives us a {\it reduced}
complex structure ${\mathcal O}_X\left/{\mathcal I}_\psi\right.$ for
$Y^\prime$.   The reduced complex subspace $\left(Y^\prime,{\mathcal
O}_X\left/{\mathcal I}_\psi\right.\right)$ is called a minimum
center of log canonical singularity.  (Usually for this technique of
minimum center of log canonical singularity one requires, in
addition, that the proper transform of $Y^\prime$ in some
appropriate blow-up $\tilde X$ of $X$ to be a nonsingular
hypersurface in $\tilde X$.)  Now one replaces $Y$ by $Y^\prime$ and
uses the vanishing of
$$
H^1\left(X,{\mathcal I}_\psi\left(mL\right)\right)=0
$$
and the surjectivity of the map
$$
\Gamma\left(X,mL\right)\to\Gamma\left(Y^\prime,\left({\mathcal
O}_X\left/{\mathcal I}_\psi\right.\right)\left(mL\right)\right)
$$
to reduce the problem to the construction of nonzero element of
$$
\Gamma\left(Y^\prime,\left({\mathcal O}_X\left/{\mathcal
I}_\psi\right.\right)\left(mL\right)\right)=\Gamma\left(Y^\prime,
mL|_{Y^\prime}\right).
$$

\bigbreak\noindent{\it Constrained Minimum Center of Log Canonical
Singularity.}  For our case at hand for the finite generation of the
canonical ring one modification has to be added in the application
of the technique of minimum center of log canonical singularity.
This modification necessitates the introduction of a new concept
which we give the name {\it constrained minimum center of log
canonical singularity} just to make it easier to refer to.  Let us
now describe the situation. In the second case of the dichotomy of
the curvature current, there is some non identically zero element
$s_0$ of $\Gamma\left(Y,mL|_Y\right)$ which is explicitly
constructed from the canonical decomposition on $Y$ of a modified
restriction of the curvature current.

\medbreak The section $s_0$ may have a nonempty zero-set $W$ in $Y$.
If we just use the technique of minimum center of log canonical
singularity without any modification, we may end up with a minimum
center of log canonical singularity $Y^\prime$ which is completely
contained inside the zero-set $W$ of $s_0$.  In such a case the
extension of $s_0|_{Y^\prime}$ is useless, because the restriction
$s_0|_{Y^\prime}$ of $s_0$ to $Y^\prime$ is identically zero on
$Y^\prime$. So we need to introduce a modification to the technique
of minimum center of log canonical singularity.  In the procedure of
using a metric $e^{-\psi}$ with least singularity to get the minimum
center of log canonical singularity, we introduce the additional
condition that the support of ${\mathcal O}_X\left/{\mathcal
I}_\psi\right.$ is not contained entirely in the zero-set $W$ of
$s_0$.  With this additional condition we can no longer require that
the structure sheaf ${\mathcal O}_X\left/{\mathcal I}_\psi\right.$
of $Y^\prime$ is reduced, but we can require that the set $E$ of
points where the structure sheaf ${\mathcal O}_X\left/{\mathcal
I}_\psi\right.$ of $Y^\prime$ fails to be reduced is entirely
contained in $W$.  We call $Y^\prime$, which is obtained from this
procedure of the additional condition, a {\it constrained minimum
center of log canonical singularity}.

\medbreak The key point about the use of a constrained minimum
center of log canonical singularity is the following.  Though the
restriction $s|_{Y^\prime}$ to $Y^\prime$ is only holomorphic on the
reduced structure of $Y^\prime$, yet since $s_0$ vanishes on $E$ we
can take a sufficiently high power $s^N$ of $s$ so that
$s^N|_{Y^\prime}$ is holomorphic on the unreduced structure
${\mathcal O}_X\left/{\mathcal I}_\psi\right.$ of $Y^\prime$.  We
now extend $s^N|_{Y^\prime}$ to $X$.  Of course, we have to replace
$m$ by $Nm$.

\bigbreak\noindent{\sc Proposition} {\it (Global Generation of the
Pluricanonical Bundle at Points of Zero Stable Vanishing Order).}
Let $X$ be a compact complex algebraic manifold of complex dimension
$n$.  Let $e^{-\varphi}$ be the metric of $K_X$ of minimum
singularity and let $\Theta_\varphi$ be its curvature current.  Then
there exist a positive integer $m_0$ such that the common zero-set
$W$ of a ${\mathbb C}$-basis of $\Gamma\left(X,m_0K_X\right)$ is
precisely the set of points where the Lelong number of
$\Theta_\varphi$ is positive.

\medbreak\noindent{\sc Proof.}  We use the technique for the Fujita
conjecture and constrained minimum centers of log canonical
singularity.  We use $W$ as the set of constraint for the
constrained minimum center of log canonical singularity.  We will
not go into further details here, because a similar but harder
situation will be handled in the proof of precisely achieving stable
vanishing order $\gamma>0$ for codimension one in the case of a
hypersurface $Y$ whose coefficient in $\Theta_\varphi$ is $\gamma$.
The only difference is that here the number $\gamma$ is replaced by
$0$ and we do the argument in the ambient space $X$ instead of in
the hypersurface $Y$ taking the place of $X$.

\medbreak We now finish the use of the temporary meaning of $X$ and
$L$ and return to the earlier meaning of $X$ and $L$.

\bigbreak\noindent{\sc Proposition}{\it (Extension of Explicitly
Constructed Section to Ambient Manifold by Constrained Minimum
Center of Log Canonical Singularity).} Let $X$ be a compact complex
algebraic manifold of complex dimension $n$. Let $e^{-\varphi}$ be
the metric of $K_X$ of minimum singularity and let $\Theta_\varphi$
be its curvature current. Let $M$ be a nonsingular hypersurface in
$X$ such that the stable vanishing order $\eta$ for $M$ is a
positive rational number. Let
$$
\left(\Theta_\varphi-\eta M\right)|_M=\sum_{j=1}\gamma_jY_j
$$
be the canonical decomposition of the closed positive
$(1,1)$-current $\left(\Theta_\varphi-\eta M\right)|_M$ on $M$ with
$J<\infty$ and each $\gamma_j$ being rational. Then the stable
vanishing order $\eta$ for $M$ is precisely achieved.

\medbreak\noindent{\sc Proof.} By the previous Proposition we find a
positive integer $m_0$ such that the common zero-set of a ${\mathbb
C}$-basis $s_1,\cdots,s_k$ of $\Gamma\left(X,m_0K_X\right)$ is
precisely the set of points where the Lelong number of
$\Theta_\varphi$ is positive. If any of the elements of
$\Gamma\left(X,m_0K_X\right)$ precisely achieves the stable
vanishing order $\eta$ for $M$, we are already done.  By replacing
$m_0$ by another sufficiently large integer, we can make the
vanishing order across $M$ of
$\sum_{j=1}^k\left|s_j\right|^{\frac{2}{m_0}}$ to be as close to
$\eta$ as prescribed (though still $>\eta$).  On the other hand, by
raising $\sum_{j=1}^k\left|s_j\right|^{\frac{2}{m_0}}$ to a positive
integral power afterwards and using interpolation, we can adjust the
stable vanishing order across $M$ to $\eta$ plus any positive
prescribed number.

\medbreak Let $L=K_X-\eta M$.  We can thus use $s_1,\cdots,s_k$ and
interpolation of metrics and their slight modifications to construct
a metric $e^{-\chi}$ of $mL-K_X=(m-1)L-(\eta+1)M$ of strictly
positive curvature current so that \begin{itemize}\item[(i)] the
zero-set of its multiplier ideal sheaf ${\mathcal I}_\chi$ is
contained in the set of points where the Lelong number of
$\Theta_\varphi$ is positive, and\item[(ii)] the generic vanishing
order of the its multiplier ideal sheaf ${\mathcal I}_\chi$ across
$M$ is precisely $1$.
\end{itemize}
We are able to fulfill Condition(i), because we can construct the
metric $e^{-\chi}$ by using the $(m-1)$-th power
$\sum_{j=1}^k\left|s_j\right|^{\frac{2}{m_0}}$ and since the common
zero-set of a ${\mathbb C}$-basis $s_1,\cdots,s_k$ of
$\Gamma\left(X,m_0K_X\right)$ is precisely the set of points where
the Lelong number of $\Theta_\varphi$ is positive.  We are able to
fulfill Condition(ii), because we have the additional order $\eta+1$
across $M$ to spare when we use the $(m-1)$-th power
$\sum_{j=1}^k\left|s_j\right|^{\frac{2}{m_0}}$ and require only the
higher generic vanishing order of $m\eta+1=(m-1)\eta+(\eta+1)$
across $M$ instead of the order $(m-1)\eta$ from the stable
vanishing order $\eta$ across $M$.

\medbreak The vanishing theorem of Kawamata-Viehweg-Nadel gives
$$
H^1\left(X,{\mathcal I}_\chi\left(mL\right)\right)=0.
$$
To get elements of $\Gamma\left(X,mL\right)$, we need to use
elements of $\Gamma\left(X,\left({\mathcal O}_X\left/{\mathcal
I}_\chi\right.\right)\left(mL\right)\right)$.   At points of $Y$ the
additional vanishing order of the multiplier ideal sheaf ${\mathcal
I}_\chi$ of $e^{-\chi}$ occurs only at points where the Lelong
number of $\Theta_\varphi$ is positive.  When we construct a
constrained minimum center $M^\prime$ of log canonical singularity
in $M$ by interpolation (with the subvariety $\cup_{j=1}^J Y_j$ as
the set of constraint), the complex structure of $M^\prime$ is
reduced outside $\cup_{j=1}^J Y_j$.  Now we can explicitly construct
the section
$$
s_0=\prod_{j=1}^k\left(s_{Y_j}\right)^{\gamma_j N}
$$
of the line bundle $NL|_M$ on the reduced structure of $M$ for some
appropriately chosen positive integer $N$, where $s_{Y_j}$ is the
canonical section of $Y_j$ on $M$.  Since $s_0$ vanishes on
$s_{Y_j}$, by replacing $N$ by a large positive integral multiple,
we can assume that the restriction of $s_0$ to $M^\prime$ can be
extended to a holomorphic section over any unreduced structure of
$M^\prime$ which is reduced outside $\cup_{j=1}^J Y_j$ and,
moreover, can be extended to an element of
$\Gamma\left(X,\left({\mathcal O}_X\left/{\mathcal
I}_\chi\right.\right)\left(NL\right)\right)$. Note that the support
of ${\mathcal O}_X\left/{\mathcal I}_\chi\right.$ may be more than
just $Y$, but its intersection with $Y$ is contained in the zero-set
$\cup_{j=1}^J Y_j$ of $s_0$ so that the last extension to an element
of $\Gamma\left(X,\left({\mathcal O}_X\left/{\mathcal
I}_\chi\right.\right)\left(NL\right)\right)$ is possible for a
sufficiently large $N$.  Since $s_0$ is nonzero at some point $P_0$
of $M^\prime$, the stable vanishing order across $M$ is precisely
achieved by $$\left(\tilde
s_0\left(s_M\right)^{N\eta}\right)^q\in\Gamma\left(X,NqK_X\right)$$
at the point $P_0$ of $M$.  Q.E.D.

\eject\centerline{\bf PART II} \centerline{\bf Illustration in Low
Dimension of the Argument of}\centerline{\bf Precise Achievement of
Stable Vanishing Order}\centerline{\bf for Higher Codimension}

\bigbreak We now illustrate the argument of the precise achievement
of the stable vanishing order for higher codimension by using the
low dimensional cases of complex surfaces and complex threefolds.
First we consider the case of surfaces.  For surfaces codimension
two means isolated points.  For isolated points for any dimension
there is a simple direct argument, which is given in the following
proposition.  After we present the case of threefolds, we will
remark on how this simple direct argument can be interpreted in the
context of the argument for any dimension.

\bigbreak\noindent{\sc Proposition} {\it(Precise Achievement of
Stable Vanishing Order at a Finite Set).}  The stable vanishing
order is automatically precisely achieved everywhere when it is
precisely achieved outside a finite set of a compact complex
algebraic manifold $X$.

\bigbreak\noindent{\sc Proof.} Let $X$ be a compact complex
algebraic manifold of general type. Let $e^{-\varphi}$ be the metric
of minimum singularity for the canonical line bundle $K_X$ of $X$.
Suppose it has been proved that the stable vanishing order is
precisely achieved outside a finite number of points
$P_1,\cdots,P_k$ of $X$ by using
$$
\sigma_1,\cdots,\sigma_\ell\in\Gamma\left(X,m_0K_X\right).
$$
We are going to show that this finite number of points must be the
empty set, otherwise there is a contradiction.  Let
$$
e^{-\psi}=\frac{1}{\sum_{j=1}^\ell\left|\sigma_j\right|^2}.
$$
We take $p\in{\mathbb N}$ sufficiently large to magnify the
discrepancy of the vanishing orders of $e^{-pm_0\varphi}$ and
$e^{-p\psi}$ so that the support of the quotient
$$
{\mathcal I}_{pm_0\varphi}\left/{\mathcal I}_{p\psi}\right.
$$
of the multiplier ideal sheaves ${\mathcal I}_{pm_0\varphi}$ and
${\mathcal I}_{p\psi}$ of $e^{-pm_0\varphi}$ and $e^{-p\psi}$
respectively is the finite set $\left\{P_1,\cdots,P_k\right\}$. We
now apply a slight modification to ${\mathcal I}_{p\psi}$ to get a
metric with strictly positive curvature current.  Let $e^{-\theta}$
be a metric of $K_X$ with strictly positive curvature current. Since
$e^{-\varphi}$ is the metric of minimum singularity for $K_X$, it
follows that when $\varepsilon>0$ is sufficiently small (which we
assume to be the case) the multiplier ideal sheaf ${\mathcal
I}_{\left(p-\varepsilon\right)\psi+\varepsilon m_0\theta}$ of the
metric $e^{-\left(\left(p-\varepsilon\right)\psi+\varepsilon
m_0\theta\right)}$ agrees with the multiplier ideal sheaves
${\mathcal I}_{pm_0\varphi}$ of $e^{-pm_0\varphi}$ on
$X-\left\{P_1,\cdots,P_k\right\}$ and the support of the quotient
$$
{\mathcal I}_{pm_0\varphi}\left/{\mathcal
I}_{\left(p-\varepsilon\right)\psi+\varepsilon m_0\theta}\right.
$$
of the multiplier ideal sheaves ${\mathcal I}_{pm_0\varphi}$ and
${\mathcal I}_{p\psi}$ of $e^{-pm_0\varphi}$ and
$e^{-\left(\left(p-\varepsilon\right)\psi+\varepsilon
m_0\theta\right)}$ respectively is the finite set
$\left\{P_1,\cdots,P_k\right\}$.  By the vanishing theorem of
Kawamata-Viehweg-Nadel, we have
$$
H^1\left(X,{\mathcal
I}_{p\psi}\left(\left(pm_0+1\right)K_X\right)\right)=0.
$$
Then the map
$$
\Gamma\left(X,{\mathcal
I}_{pm_0\varphi}\left(\left(pm_0+1\right)K_X\right)\right)\to
\Gamma\left(X,\left({\mathcal I}_{pm_0\varphi}\left/{\mathcal
I}_{\left(p-\varepsilon\right)\psi+\varepsilon
m_0\theta}\right.\right)\left(\left(pm_0+1\right)K_X\right)\right)
$$
is surjective.  Note that for this we do not need the vanishing of
the cohomology group
$$
H^1\left(X,{\mathcal
I}_{pm_0\varphi}\left(\left(pm_0+1\right)K_X\right)\right).
$$
Since
$$
\Gamma\left(X,\left({\mathcal I}_{pm_0\varphi}\left/{\mathcal
I}_{\left(p-\varepsilon\right)\psi+\varepsilon
m_0\theta}\right.\right)\left(\left(pm_0+1\right)K_X\right)\right)\approx
\bigoplus_{j=1}^k\left({\mathcal I}_{pm_0\varphi}\left/{\mathcal
I}_{\left(p-\varepsilon\right)\psi+\varepsilon
m_0\theta}\right.\right)_{P_j},
$$
it follows from Nakayama's lemma and the surjectivity of
$$
\Gamma\left(X,{\mathcal
I}_{pm_0\varphi}\left(\left(pm_0+1\right)K_X\right)\right)\to
\bigoplus_{j=1}^k\left({\mathcal I}_{pm_0\varphi}\left/{\mathcal
I}_{\left(p-\varepsilon\right)\psi+\varepsilon
m_0\theta}\right.\right)_{P_j},
$$
that the map
$$
\Gamma\left(X,{\mathcal
I}_{pm_0\varphi}\left(\left(pm_0+1\right)K_X\right)\right)\to
\bigoplus_{j=1}^k\left({\mathcal I}_{pm_0\varphi}\right)_{P_j}
$$
is surjective.  This actually gives a contradiction, because the
stable vanishing order of $\left(pm_0+1\right)K_X$ should be given
by $e^{-\left(pm_0+1\right)\varphi}$ instead of by
$e^{-pm_0\varphi}$.  Q.E.D.

\bigbreak\noindent{\it Remarks.} (a) With this proposition, to get
the finite generation of the canonical ring for the case of a
compact complex algebraic surface of general type it suffices to
show that the stable vanishing order is precisely achieved at
codimension one.

\medbreak\noindent(b) For the analytic proof of the finite
generation of the canonical ring, when we get down to the point of
having already verified the precise achievement of the stable
vanishing order outside a finite set of points, we do not need to
blow up the points to reduce the argument to the case of a
hypersurface in the new blown-up manifold.  Therefore the difficulty
does not exist, in the case of a surface, of blowing up a point
$P_0$ to get a curve $C_1$ and then locating some bad point $P_1$
(where the precise achievement of stable vanishing order fails) in
the curve $C_1$ and then blowing up $P_1$ to get a curve $C_2$ and
locating some bad point $P_2$ (where the precise achievement of
stable vanishing order fails) in the curve $C_2$ and then blowing up
$P_2$ and possibly finally ending up with an unending infinite
sequence of bad points, each in a tower of successively blown-up
surfaces.

\bigbreak\noindent{\it Higher Codimension Argument for Threefold
Case.} Let $X$ be a complex complex algebraic threefold of general
type. Let $e^{-\varphi}$ be the metric of minimum singularity for
the canonical line bundle $K_X$ of $X$.  Suppose it has been proved
that the stable vanishing order is precisely achieved outside a
curve $C=\bigcup_jC_j$ (where each $C_j$ is irreducible) by
$$
s_1,\cdots,s_k\in\Gamma\left(X,m_0K_X\right).
$$
Note that here we have used the Proposition given above to rule out
the possibility that, besides at the curve $C$, the stable vanishing
order may fail to be precisely achieved at a finite set of points in
$X-C$.

\medbreak On $X$ we have the canonical decomposition of
$$\Theta_{K_X}=\sum_{j=1}^J\gamma_j Y_j+R$$
of the curvature current of the metric of minimum singularity for
$K_X$.  We are skipping the diophantine argument which is explained
in detail in [Siu 2006] and consider as verified the rationality of
each $\gamma_j$. By replacing $m_0$ by an appropriate integral
multiple, we can assume without loss of generality that each
$m_0\gamma_j$ is a positive integer.  Let
$s^*=\prod_{j=1}^J\left(s_{Y_j}\right)^{m_0\gamma_j}$, where
$s_{Y_j}$ is the canonical section of $Y_j$.  Let
$L=K_X-\sum_{j=1}^J\gamma_j Y_j$.  Since we can replace
$\sigma_1,\cdots,\sigma_\ell$ by
$$
\frac{s_1}{s^*},\cdots,\frac{s_k}{s^*}\in\Gamma\left(X,m_0L\right)
$$
and the essence of the rest of the argument that is to follows does
not change, for notational simplicity we are going to assume that
all $\gamma_j=0$ so that $L=K_X$.  Also, for notational simplicity
we are going to assume that the curve $C$ is irreducible.

\medbreak \medbreak For
$\sigma=\left(\sigma_1,\cdots,\sigma_k\right)\in{\mathbb
C}^k-\left\{0\right\}$ let $s_\sigma=\sum_{j=1}^k\sigma_j s_j$ and
let $S_\sigma$ be the surface in $X$ defined by $s_\sigma$.  We
consider those $\sigma$ for which $S_\sigma$ is irreducible across
which $s_\sigma$ vanishes to order $1$.  By considering the blow-up
of $X$ by the ideal generated by $s_1,\cdots,s_k$ and the precise
achievement of stable vanishing order for codimension one, after
replacing $m_0$ by a positive integral multiple if necessary, we
have the following situation.
\begin{itemize}\item[(a)] For $\sigma=\left(\sigma_1,\cdots,\sigma_k\right)\in{\mathbb
C}^k-\left\{0\right\}$ for which $S_\sigma$ is irreducible where
$s_\sigma$ vanishes to order $1$, there exist
\begin{itemize}\item[(i)] some $\tau=\left(\tau_1,\cdots,\tau_k\right)\in{\mathbb
C}^k-\left\{0\right\}$ and\item[(ii)] some finite subset
$Z_{\tau,\sigma}$ of $C$\end{itemize} such that the section $s_\tau$
is not identically zero on $S_\sigma$ and the multi-valued section
$$
\frac{s_\tau|_{S_\sigma}}{\left(s_{C,\sigma}\right)^{\eta_\sigma}}
$$
on a neighborhood of $C-Z_{\tau,\sigma}$ in $S_\sigma$ is nonzero at
points of $C-Z_{\tau,\sigma}$, where $\eta_\sigma$ is the stable
vanishing order at $C$ on the surface $S_\sigma$ and $s_{C,\sigma}$
is the canonical section of $C$ in $S_\sigma$.
\item[(b)] One has the second case of the dichotomy for the
curvature current
$$
\left.\left(\Theta_{K_X}|_{S_\sigma}-\eta_\sigma\left[C\right]\right)\right|_C
$$
on $C$.
\item[(c)] The multi-valued section
$$
\frac{s_\tau|_{S_\sigma}}{\left(s_{C,\sigma}\right)^{\eta_\sigma}}
$$
on $C$ is constructed from the second case of the dichotomy for the
curvature current
$$
\left.\left(\Theta_{K_X}|_{S_\sigma}-\eta_\sigma\left[C\right]\right)\right|_C
$$
on $C$.

\item[(c)] The stable vanishing order $\eta_\sigma$ for $C$ on $S_\sigma$ is
achieved by $s_\tau$ at points of $C-Z_{\tau,\sigma}$.
\end{itemize}
Note that, because of Condition(b) and Condition(c) we can choose
$Z_{\tau,\sigma}$ independent of $\tau$ so that
$Z_{\tau,\sigma}=Z_\sigma$ for some finite subset $Z_\sigma$ of $C$
depending only on $\sigma$.  To finish the proof of the precise
achievement of the stable vanishing order for the threefold $X$, it
suffices to show that there is a finite subset $Z$ of $C$ such that
every $Z_\sigma$ is contained in $Z$.

\medbreak At a regular point $P_0$ of $C$, there is some open
neighborhood $U$ of $P_0$ in $X$ such that
\begin{itemize}\item[(i)] the pair $(U,C)$ is biholomorphic to the pair
$\left(\Delta^3,\left\{(0,0)\right\}\times\Delta\right)$, where
$\Delta$ is the open unit disk in ${\mathbb C}$, and
\item[(ii)] $s_j|U$ is represented by a holomorphic function
$f_j\left(z_1,z_2,t\right)$ on $\Delta^3$ for $1\leq j\leq k$.
\end{itemize}
The existence of the finite subset $Z$ of $C$ with the property that
$Z_{\tau,\sigma}\subset Z$ follows from Lemma 2 given below.

\bigbreak\noindent{\it Lemma 1.}  Let $f_0,\cdots,f_k$ be
holomorphic function germs on ${\mathbb C}^2$ at the origin so that
the origin is the common zero-set of any two of the holomorphic
function germs $f_0,\cdots,f_k$.  Let $C_0$ be the complex curve
germ at the origin defined by $f_0=0$, which is assumed to be
irreducible and across which $f_0$ vanishes to order $1$.  Then the
following numbers are the same.
\begin{itemize}
\item[(i)] The multiplicity of the ideal $\sum_{j=0}^k{\mathcal
O}_{{\mathbb C}^2,0}f_j$ at the origin.
\item[(ii)] The dimension over ${\mathbb C}$ of
$$
{\mathcal O}_{{\mathbb C}^2,0}\left/\sum_{j=0}^k{\mathcal
O}_{{\mathbb C}^2,0}f_j\right..
$$
\item[(iii)] The Lelong number $\lambda$ of
$$
\hat\Theta:=\left.\left(\frac{\sqrt{-1}}{2\pi}\partial\bar\partial\log\sum_{j=1}^k\left|f_j\right|^2\right)\right|_{C_0}
$$
on $C_0$ at the origin.  Here the Lelong number $\lambda$ means the
Lelong number of the pullback of the closed positive $(1,1)$-current
$\hat\Theta$ to the normalization of $C_0$.

\item[(iv)] The number $\eta$ such that
$$
\frac{\sum_{j=1}^k\left|f_j\right|^2}{\left(\sum_{j=1}^2\left|z_j\right|^2\right)^\eta}
$$
is bounded between two positive numbers near the origin on $C_0$.
\end{itemize}

\bigbreak\noindent{\it Lemma 2.}  Let $k\geq 2$ and let
$f_j(z_1,z_2,t)$, for $1\leq j\leq k$, be holomorphic functions on
the tri-disk $\Delta^3$ with coordinates $z_1,z_2,t$. Assume that
the common zero-set of any two of $f_1,\cdots,f_k$ is
$\left\{(0,0)\right\}\times\Delta=\left\{z_1=z_2=0\right\}$. For any
$k$-tuple $\left(a_1,\cdots,a_k\right)$ of complex numbers not all
zero let $S_{a_1,\cdots,a_k}$ be the zero-set of $\sum_{j=1}^k
a_jf_j$.  Then there exists a discrete subset $Z$ of
$\left\{(0,0)\right\}\times\Delta$ with the following property. For
any $k$-tuple $\left(a_1,\cdots,a_k\right)$ of complex numbers not
all zero there exists some nonnegative number
$\gamma_{a_1,\cdots,a_k}$ such that
$$
\frac{\sum_{j=1}^k\left|f_j\right|^2}{\left(\sum_{j=1}^2\left|z_j\right|^2\right)^{\gamma_{a_1,\cdots,a_k}}}$$
is continuous nonzero on some neighborhood of
$\left(S_{a_1,\cdots,a_k}\cap\left\{(0,0)\right\}\times\Delta\right)-Z$
in $S_{a_1,\cdots,a_k}$ if $S_{a_1,\cdots,a_k}$ is irreducible and
$\sum_{j=1}^k a_jf_j$ vanishes to order $1$ across
$S_{a_1,\cdots,a_k}$.

\medbreak\noindent{\it Proof.}  Observe that for fixed $t$ if
$a_1\not=0$ then each of $\sum_{j=1}^k\left|f_j\right|^2$ and
$\left|a_1f_1+\cdots+a_kf_k\right|^2+\sum_{j=2}^k\left|f_j\right|^2$
is bounded by a positive constant times the other on some
neighborhood of the origin in ${\mathbb C}^2$.  Use Lemma 1. Q.E.D.

\bigbreak\noindent{\it Remarks.}  (a) We would like to highlight the
intuitive geometric reason for the existence of a discrete set $Z$
in $C$ such that the ``bad set'' $Z_\sigma$ in $C$ for the
restriction of the closed positive $(1,1)$-current $\Theta_{K_X}$
restricted to $S_\sigma$ is contained in $Z$.   The key point is
that the surface $S_\sigma$ is sliced out by a ${\mathbb C}$-linear
combination of the pluricanonical sections $s_1,\cdots,s_k$ and that
these pluricanonical sections $s_1,\cdots,s_k$ have the property
that the ``bad set'' $Z_\sigma$ of $C$ can be described by the extra
vanishing of $\sum_{j=1}^k\left|s_j\right|^2$ beyond their generic
vanishing order at points of $C$.

\medbreak\noindent(b) We need to restrict $\Theta_{K_X}$ to
$S_\sigma$, because we have to take away the vanishing order of
$\Theta_{K_X}$ at $C$ and it is only for codimension one we can take
away the vanishing order.  The vanishing order on different
$S_\sigma$ may be different.

\medbreak\noindent(c) The ``bad set'' $Z_\sigma$ in $C$ describes
the points where the relative position of the pair of two Artinian
subschemes in the normal direction of $C$ jumps. One Artinian
subscheme comes from the restriction of $s_1,\cdots,s_k$ to a local
surface $T$ normal to $C$ and the other one comes from
$s_1,\cdots,s_k$ plus a ${\mathbb C}$-linear combination $\sigma_1
s_1+\cdots+\sigma_k s_k$ after restricting them to $T$. The main
point is that there is a finite subset $Z$ of $C$ such that
$Z_\sigma\subset Z$ with $Z$ independent of the choice of $\sigma$.

\bigbreak\noindent{\sc References.}

\bigbreak\noindent[Angehrn-Siu 1995] U. Angehrn and Y.-T. Siu,
Effective freeness and point separation for adjoint bundles. {\it
Invent. Math.} {\bf 122} (1995), 291--308.

\medbreak\noindent[Ohsawa 2002] T. Ohsawa, A precise $L^2$ division
theorem. In: {\it Complex Geometry} (Selected papers dedicated to
Hans Grauert from the International Conference on Analytic and
Algebraic Methods in Complex Geometry held in G\"ottingen, April
3--8, 2000), {\it ed.} I. Bauer, F, Catanese, Y. Kawamata, Th.
Peternell and Y.-T. Siu. Springer-Verlag, Berlin, 2002, pp.185--191.

\medbreak\noindent [Siu 1998] Y.-T. Siu, Invariance of plurigenera.
Invent. Math. {\bf 134} (1998), no. 3, 661--673.

\medbreak\noindent [Siu 2002] Y.-T. Siu, Extension of twisted
pluricanonical sections with plurisubharmonic weight and invariance
of semipositively twisted plurigenera for manifolds not necessarily
of general type. In: {\it Complex Geometry} (Selected papers
dedicated to Hans Grauert from the International Conference on
Analytic and Algebraic Methods in Complex Geometry held in
G\"ottingen, April 3--8, 2000), {\it ed.} I. Bauer, F, Catanese, Y.
Kawamata, Th. Peternell and Y.-T. Siu. Springer-Verlag, Berlin,
2002, pp.223--277.

\medbreak\noindent [Siu 2004] Y.-T. Siu, Invariance of plurigenera
and torsion-freeness of direct image sheaves of pluricanonical
bundles. In: {\it Finite or infinite dimensional complex analysis
and applications, Adv. Complex Anal. Appl.}, {\bf 2}, Kluwer Acad.
Publ., Dordrecht, 2004, pp. 45--83.

\medbreak\noindent [Siu 2006] Y.-T. Siu, A general non-vanishing
theorem and an analytic proof of the finite generation of the
canonical ring, arXiv:math/0610740.

\medbreak\noindent[Shokurov 1985] V.~V. Shokurov, A nonvanishing
theorem. {\it Izv. Akad. Nauk SSSR Ser. Mat.} {\bf 49} (1985),
635--651.

\medbreak\noindent[Skoda 1972] H. Skoda, Application des techniques
$L^2$ \`a la th\'eorie des id\'eaux d'une alg\`ebre de fonctions
holomorphes avec poids. {\it Ann. Sci. \'Ecole Norm. Sup.} {\bf 5}
(1972), 545-579.

\end{document}